\documentclass[11pt,a4paper]{article}

\usepackage{geometry}
\usepackage[T1]{fontenc}
\usepackage[utf8]{inputenc}
\usepackage{authblk}

\usepackage{graphics}
\usepackage{tikz,pgf}
\usetikzlibrary{decorations.markings, decorations.pathmorphing}
\usetikzlibrary{arrows,decorations,backgrounds,scopes,plotmarks}
\usepackage{subfigure}
\usepackage{here}
\usepackage{caption}
\usepackage{float}

\usepackage[english]{babel}
\usepackage[all]{xy}
\usepackage{marvosym}
\usepackage{stmaryrd,mathrsfs,amsmath,amssymb,bm,amsfonts}
\usepackage{enumerate}
\usepackage{amsmath,amsthm,amssymb,amscd}

\newtheorem{thm}{Theorem}[section]
\newtheorem{defn}[thm]{Definition}
\newtheorem{prop}[thm]{Proposition}
\newtheorem{lem}[thm]{Lemma}

\begin{document}

\title{Remark on topological nature of upward planarity}

\author{Xuexing Lu}

\affil{\small School of Mathematical Sciences, Zaozhuang University, Zaozhuang, China}

\renewcommand\Authands{ and }
\maketitle

\begin{abstract}
The notion of an upward plane graph in graph theory and that of a progressive plane graph (or plane string diagram) in category theory are essentially the same thing.
In this paper, we combine the ideas in graph theory and category theory to explain why and in what sense upward planarity is a topological property. The main result is that two upward planar drawings of an acyclic directed graph are equivalent (connected by a deformation) if and only if they are connected by a planar isotopy which preserves the orientation and polarization of $G$. This result gives a positive answer to Selinger's conjectue, whose strategy is different from the solution recently given by Delpeuch and Vicary.
\end{abstract}

\text{\textit{Keywords}:  upward planar drawing, progressive plane graph, recumbent isotopy}\\

%%%%%%%%%%%%%%%%%%%%%%%%%%%%%%%%%%%%%%%%%%%%%%%%%%%%%%%%%%%%%%%%%
%\tableofcontents
%%%%%%%%%%%%%%%%%%%%%%%%%%%%%%%%%%%%%%%%%%%%%%%%%%%%%%%%%%%%%%%

\section{Motivation}
It is well-known that two planar drawings on a surface are called equivalent if they are connected by a planar isotopy on the surface, which characterizes planarity as a purely topological property of graphs.
Just as planarity means to graphs, it is believed that upward planarity is a fundament topological property of acyclic directed graphs and orders. If this is true, then we should answer such two questions that $(Q1)$ what is the proper definition of equivalence relation of upward planar drawings, which  characterizes upward planarity as a purely topological property; and $(Q2)$ what is the proper definition of an upward planar drawing on a surface of higher genus.    Admirably, the answers to $(Q1)$ and $(Q2)$ are already implicit in the famous work of Joyal and Street \cite{[JS91]} on string diagrams.

In this paper, we will only focus on $(Q1)$ and reinterpret the first part of Joyal-Street's work --- the graphical calculus for monoidal categories --- as a proper topological theory of upward planar drawings. The other two parts of their work --- the graphical calculus for symmetric monoidal categories and the graphical calculus for braided monoidal categories --- should be related to $(Q2)$ and the theory of virtual knots, respectively. We will elaborate on these connections elsewhere.

\section{Upward plane graph and progressive plane graph}

\begin{figure}[H]
\centering
\begin{tikzpicture}[scale=0.25]

\node (v2) at (0.5,2) {};
\node (v1) at (-3,-3) {};
\node (v3) at (3,-2.5) {};
\node (v4) at (7,0) {};
\node (v8) at (0.5,-2) {};
\node (v5) at (-0.5,-7.5) {};
\node (v6) at (3.5,-5.5) {};
\node (v7) at (6,-8.5) {};
\draw  (-3,-3) -- (0.5,2)[postaction={decorate, decoration={markings,mark=at position .5 with {\arrowreversed[black]{stealth}}}}];
\draw  (3,-2.5) -- (0.5,2)[postaction={decorate, decoration={markings,mark=at position .5 with {\arrowreversed[black]{stealth}}}}];
\draw  (3,-2.5) -- (7,0)[postaction={decorate, decoration={markings,mark=at position .5 with {\arrowreversed[black]{stealth}}}}];
\draw  (-0.5,-7.5) -- (3.5,-5.5)[postaction={decorate, decoration={markings,mark=at position .5 with {\arrowreversed[black]{stealth}}}}];
\draw   (6,-8.5) -- (3.5,-5.5)[postaction={decorate, decoration={markings,mark=at position .5 with {\arrowreversed[black]{stealth}}}}];
\draw  (-0.5,-7.5) -- (-3,-3)[postaction={decorate, decoration={markings,mark=at position .5 with {\arrowreversed[black]{stealth}}}}];
\draw  (-0.5,-7.5) -- (0.5,-2)[postaction={decorate, decoration={markings,mark=at position .5 with {\arrowreversed[black]{stealth}}}}];
\draw  (3.5,-5.5)-- (0.5,-2)[postaction={decorate, decoration={markings,mark=at position .5 with {\arrowreversed[black]{stealth}}}}];
\draw   (6,-8.5) -- (7,0)[postaction={decorate, decoration={markings,mark=at position .5 with {\arrowreversed[black]{stealth}}}}];
\node (v9) at (2,-11) {};
\draw  (2,-11) --  (6,-8.5)[postaction={decorate, decoration={markings,mark=at position .5 with {\arrowreversed[black]{stealth}}}}];
\draw  (2,-11) -- (-0.5,-7.5)[postaction={decorate, decoration={markings,mark=at position .5 with {\arrowreversed[black]{stealth}}}}];
\draw  (3.5,-5.5) -- (7,0)[postaction={decorate, decoration={markings,mark=at position .5 with {\arrowreversed[black]{stealth}}}}];
\draw[fill] (v1) circle [radius=0.2];
\draw[fill] (v2) circle [radius=0.2];
\draw[fill] (v3) circle [radius=0.2];
\draw[fill] (v4) circle [radius=0.2];
\draw[fill] (v5) circle [radius=0.2];
\draw[fill] (v6) circle [radius=0.2];
\draw[fill] (v7) circle [radius=0.2];
\draw[fill] (v8) circle [radius=0.2];
\draw[fill] (v9) circle [radius=0.2];
\draw  plot[smooth, tension=.7] coordinates {(7,0) (8,-2.5) (8.5,-5.5) (7.5,-7.5) (6,-8.5)}[postaction={decorate, decoration={markings,mark=at position .5 with {\arrow[black]{stealth}}}}];
\end{tikzpicture}
\caption{An upward planar drawing or upward plane graph}
\label{upg}
\end{figure}
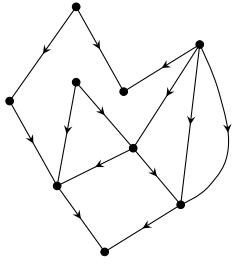

A planar drawing of a directed graph is called \textbf{upward} if all its edges monotonically decrease in the vertical direction (or other fixed direction), see Fig \ref{upg} for an example. Clearly, a necessary condition for a directed graph to have an upward planar drawing is that it is acyclic. An acyclic directed graph together with an upward planar drawing is called an \textbf{upward plane graph}. An acyclic directed graph is called an \textbf{upward planar graph} if it admits an upward planar drawing. Upward plane graphs have been extensively studied in the fields of graph theory, graph drawing algorithms, and ordered set theory (see, e.g., \cite{[GT95]} for a review).

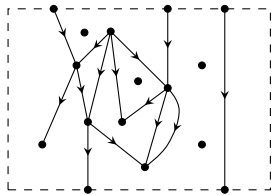
\begin{figure}[H]
\centering
\begin{tikzpicture}[scale=0.3]
\draw [dashed] (-3,1.5) rectangle (8.5,-6.5);
\node [above](v1) at (-1,1.5) {};
\node (v2) at (0,-1) {};
\node (v6) at (-1.5,-4.5) {};
\node [left] at (-1.5,-4.5) {};
\node (v3) at (1.5,0.5) {};
\node (v8) at (4,-2) {};
\node (v4) at (0.5,-3.5) {};
\node [below](v5) at (0.5,-6.5) {};
\node [above](v7) at (4,1.5) {};
\draw  (-1,1.5)  -- (0,-1)[postaction={decorate, decoration={markings,mark=at position .5 with {\arrow[black]{stealth}}}}];
\draw   (1.5,0.5) -- (0,-1)[postaction={decorate, decoration={markings,mark=at position .5 with {\arrow[black]{stealth}}}}];
\draw   (1.5,0.5)--  (0.5,-3.5)[postaction={decorate, decoration={markings,mark=at position .5 with {\arrow[black]{stealth}}}}];
\draw  (0.5,-3.5) -- (v5)[postaction={decorate, decoration={markings,mark=at position .5 with {\arrow[black]{stealth}}}}];
\draw  (0,-1) -- (-1.5,-4.5)[postaction={decorate, decoration={markings,mark=at position .5 with {\arrow[black]{stealth}}}}];
\draw  (0,-1) -- (0.5,-3.5)[postaction={decorate, decoration={markings,mark=at position .5 with {\arrow[black]{stealth}}}}];
\draw  (4,1.5)-- (4,-2)[postaction={decorate, decoration={markings,mark=at position .5 with {\arrow[black]{stealth}}}}];
\draw   (1.5,0.5)-- (4,-2)[postaction={decorate, decoration={markings,mark=at position .5 with {\arrow[black]{stealth}}}}];
\node (v9) at (3,-5.5) {};
\draw  (0.5,-3.5) -- (3,-5.5)[postaction={decorate, decoration={markings,mark=at position .5 with {\arrow[black]{stealth}}}}];
\draw  (4,-2) -- (3,-5.5)[postaction={decorate, decoration={markings,mark=at position .5 with {\arrow[black]{stealth}}}}];
\node (v10) at (2,-3.5) {};
\draw   (1.5,0.5) -- (2,-3.5)[postaction={decorate, decoration={markings,mark=at position .5 with {\arrow[black]{stealth}}}}];
\draw  (4,-2) -- (2,-3.5)[postaction={decorate, decoration={markings,mark=at position .5 with {\arrow[black]{stealth}}}}];

\draw [fill](v2) circle [radius=0.15];
\draw [fill](v3) circle [radius=0.15];
\draw [fill](v4) circle [radius=0.15];
\draw [fill](v6) circle [radius=0.15];
\draw [fill](v8) circle [radius=0.15];
\draw [fill](v9) circle [radius=0.15];
\draw [fill](v10) circle [radius=0.15];

\node at (-1,0.5) {};
\node at (4.5,0) {};
\node at (0,-5.5) {};
\node at (-1.5,-3) {};

\draw(6.5,1.5)--(6.5,-6.5)[postaction={decorate, decoration={markings,mark=at position .5 with {\arrow[black]{stealth}}}}];
\node [above]at(6.5,1.5) {};
\node [below]at (6.5,-6.5) {};
\node at (7,-3) {};

\draw [fill] (-1,1.5) circle [radius=0.15];
\draw [fill](0.5,-6.5) circle [radius=0.15];
\draw [fill](4,1.5) circle [radius=0.15];
\draw [fill](6.5,-6.5) circle [radius=0.15];
\draw [fill](6.5,1.5) circle [radius=0.15];
\draw [fill](5.5,-1) circle [radius=0.15];
\draw [fill](2.7,-1.7) circle [radius=0.15];
\draw [fill](0.35,0.45) circle [radius=0.15];
\draw  plot[smooth, tension=.7] coordinates {(4,-2) (4.5,-3) (4,-4.5) (3,-5.5)}[postaction={decorate, decoration={markings,mark=at position .5 with {\arrow[black]{stealth}}}}];
\node at (5.5,-4.5) {};
\draw [fill](5.5,-4.5) circle [radius=0.15];
\end{tikzpicture}
\caption{A progressive plane graph or plane string diagram}
\label{psd}
\end{figure}

A \textbf{progressive plane graph}, or \textbf{plane string diagram} is an upward planar drawing of an acyclic directed graph (possibly with isolated vertices) in a plane box such that
$(1)$ all vertices drawn on the horizontal boundaries  are leaves (vertices of degree one);
$(2)$ no vertices are drawn on the vertical boundaries, see Fig \ref{psd} for an example.   It was introduced by Joyal and Street \cite{[JS88],[JS91]} as a core notion in their graphical calculus for monoidal categories, which provides a firm foundation for many explorations of graphical notations in mathematics and physics.

It is not difficult to see that the notion of an upward plane graph and that of a progressive plane graph are essentially the same thing \cite{[HLY16]}.
A progressive plane graph is just an upward plane graph satisfying the so-called \textbf{boxed condition} which means that $(1)$ the graph is drawn in a plane box;
$(2)$ all vertices drawn on the horizontal boundaries are leaves; and $(3)$ no vertices are drawn on the vertical boundaries.

\begin{figure}[H]
\centering
$$
\begin{matrix}\begin{matrix}
\begin{tikzpicture}[scale=0.4]

\draw [loosely dashed] (-3,1.5) rectangle (8.5,-6.5);
\node [above](v1) at (-1,1.5) {};
\node (v2) at (0,-1) {};
\node (v6) at (-1.5,-4.5) {};
\node [left] at (-1.5,-4.5) {};
\node (v3) at (1.5,0.5) {};
\node (v8) at (4,-2) {};
\node (v4) at (0.5,-3.5) {};
\node [below](v5) at (0.5,-6.5) {};
\node [above](v7) at (4,1.5) {};
\draw  (-1,1.5)  -- (0,-1)[postaction={decorate, decoration={markings,mark=at position .5 with {\arrow[black]{stealth}}}}];
\draw   (1.5,0.5) -- (0,-1)[postaction={decorate, decoration={markings,mark=at position .5 with {\arrow[black]{stealth}}}}];
\draw   (1.5,0.5)--  (0.5,-3.5)[postaction={decorate, decoration={markings,mark=at position .5 with {\arrow[black]{stealth}}}}];
\draw  (0.5,-3.5) -- (v5)[postaction={decorate, decoration={markings,mark=at position .5 with {\arrow[black]{stealth}}}}];
\draw  (0,-1) -- (-1.5,-4.5)[postaction={decorate, decoration={markings,mark=at position .5 with {\arrow[black]{stealth}}}}];
\draw  (0,-1) -- (0.5,-3.5)[postaction={decorate, decoration={markings,mark=at position .5 with {\arrow[black]{stealth}}}}];
\draw  (4,1.5)-- (4,-2)[postaction={decorate, decoration={markings,mark=at position .5 with {\arrow[black]{stealth}}}}];
\draw   (1.5,0.5)-- (4,-2)[postaction={decorate, decoration={markings,mark=at position .5 with {\arrow[black]{stealth}}}}];
\node (v9) at (3,-5.5) {};
\draw  (0.5,-3.5) -- (3,-5.5)[postaction={decorate, decoration={markings,mark=at position .5 with {\arrow[black]{stealth}}}}];
\draw  (4,-2) -- (3,-5.5)[postaction={decorate, decoration={markings,mark=at position .5 with {\arrow[black]{stealth}}}}];
\node (v10) at (2,-3.5) {};
\draw   (1.5,0.5) -- (2,-3.5)[postaction={decorate, decoration={markings,mark=at position .5 with {\arrow[black]{stealth}}}}];
\draw  (4,-2) -- (2,-3.5)[postaction={decorate, decoration={markings,mark=at position .5 with {\arrow[black]{stealth}}}}];

\draw [fill](v2) circle [radius=0.11];
\draw [fill](v3) circle [radius=0.11];
\draw [fill](v4) circle [radius=0.11];
\draw [fill](v6) circle [radius=0.11];
\draw [fill](v8) circle [radius=0.11];
\draw [fill](v9) circle [radius=0.11];
\draw [fill](v10) circle [radius=0.11];

\node at (-1,0.5) {};
\node at (4.5,0) {};
\node at (0,-5.5) {};
\node at (-1.5,-3) {};

\draw(6.5,1.5)--(6.5,-6.5)[postaction={decorate, decoration={markings,mark=at position .5 with {\arrow[black]{stealth}}}}];
\node [above]at(6.5,1.5) {};
\node [below]at (6.5,-6.5) {};
\node at (7,-3) {};

\draw [fill] (-1,1.5) circle [radius=0.11];
\draw [fill](0.5,-6.5) circle [radius=0.11];
\draw [fill](4,1.5) circle [radius=0.11];
\draw [fill](6.5,-6.5) circle [radius=0.11];
\draw [fill](6.5,1.5) circle [radius=0.11];

\draw [fill](5.5,-1) circle [radius=0.11];
\draw [fill](2.7,-1.7) circle [radius=0.11];

\draw [fill](0.35,0.45) circle [radius=0.11];

\draw  plot[smooth, tension=.7] coordinates {(4,-2) (4.5,-3) (4,-4.5) (3,-5.5)}[postaction={decorate, decoration={markings,mark=at position .5 with {\arrow[black]{stealth}}}}];
\node at (5.5,-4.5) {};
\draw [fill](5.5,-4.5) circle [radius=0.11];
\end{tikzpicture}
\end{matrix}&&\begin{matrix}\Longrightarrow\end{matrix}&&\begin{matrix}
\begin{tikzpicture}[scale=0.4]

\draw [loosely dashed] (-3,1.5) rectangle (8.5,-6.5);
\node [above](v1) at (-1,1.5) {};
\node (v2) at (0,-1) {};
\node (v6) at (-1.5,-4.5) {};
\node [left] at (-1.5,-4.5) {};
\node (v3) at (1.5,0.5) {};
\node (v8) at (4,-2) {};
\node (v4) at (0.5,-3.5) {};
\node [below](v5) at (0.5,-6.5) {};
\node [above](v7) at (4,1.5) {};
\draw  (-1,1.5)  -- (0,-1)[postaction={decorate, decoration={markings,mark=at position .5 with {\arrow[black]{stealth}}}}];
\draw   (1.5,0.5) -- (0,-1)[postaction={decorate, decoration={markings,mark=at position .5 with {\arrow[black]{stealth}}}}];
\draw   (1.5,0.5)--  (0.5,-3.5)[postaction={decorate, decoration={markings,mark=at position .5 with {\arrow[black]{stealth}}}}];
\draw  (0.5,-3.5) -- (v5)[postaction={decorate, decoration={markings,mark=at position .5 with {\arrow[black]{stealth}}}}];
\draw  (0,-1) -- (-1.5,-4.5) node (v19) {}[postaction={decorate, decoration={markings,mark=at position .5 with {\arrow[black]{stealth}}}}];
\draw  (0,-1) node (v16) {} -- (0.5,-3.5)[postaction={decorate, decoration={markings,mark=at position .5 with {\arrow[black]{stealth}}}}];
\draw  (4,1.5)-- (4,-2)[postaction={decorate, decoration={markings,mark=at position .5 with {\arrow[black]{stealth}}}}];
\draw   (1.5,0.5)-- (4,-2)[postaction={decorate, decoration={markings,mark=at position .5 with {\arrow[black]{stealth}}}}];
\node (v9) at (3,-5.5) {};
\draw  (0.5,-3.5) -- (3,-5.5)[postaction={decorate, decoration={markings,mark=at position .5 with {\arrow[black]{stealth}}}}];
\draw  (4,-2) -- (3,-5.5) node (v22) {}[postaction={decorate, decoration={markings,mark=at position .5 with {\arrow[black]{stealth}}}}];
\node (v10) at (2,-3.5) {};
\draw   (1.5,0.5) node (v18) {} -- (2,-3.5)[postaction={decorate, decoration={markings,mark=at position .5 with {\arrow[black]{stealth}}}}];
\draw  (4,-2) -- (2,-3.5) node (v21) {}[postaction={decorate, decoration={markings,mark=at position .5 with {\arrow[black]{stealth}}}}];

\draw [fill](v2) circle [radius=0.11];
\draw [fill](v3) circle [radius=0.11];
\draw [fill](v4) circle [radius=0.11];
\draw [fill](v6) circle [radius=0.11];
\draw [fill](v8) circle [radius=0.11];
\draw [fill](v9) circle [radius=0.11];
\draw [fill](v10) circle [radius=0.11];

\node at (-1,0.5) {};
\node at (4.5,0) {};
\node at (0,-5.5) {};
\node at (-1.5,-3) {};

\draw(6.5,1.5)--(6.5,-6.5)[postaction={decorate, decoration={markings,mark=at position .5 with {\arrow[black]{stealth}}}}];
\node [above]at(6.5,1.5) {};
\node [below]at (6.5,-6.5) {};
\node at (7,-3) {};

\draw [fill] (-1,1.5) circle [radius=0.11];

\draw [fill](0.5,-6.5) circle [radius=0.11];
\draw [fill](4,1.5) circle [radius=0.11];
\draw [fill](6.5,-6.5) circle [radius=0.11];
\draw [fill](6.5,1.5) circle [radius=0.11];

\draw [fill](5.5,-1) circle [radius=0.11];
\draw [fill](5.5,0) circle [radius=0.11];
\draw [densely dotted]  (5.5,0)-- (5.5,-1)[postaction={decorate, decoration={markings,mark=at position .55 with {\arrow[black]{stealth}}}}];

\draw [fill](2.4,-1.2) circle [radius=0.11];
\draw [fill](2.4,-2.2) circle [radius=0.11];
\draw [densely dotted]  (2.4,-1.2)-- (2.4,-2.2)[postaction={decorate, decoration={markings,mark=at position .55 with {\arrow[black]{stealth}}}}];

\draw [fill](0.15,0.9) circle [radius=0.11];
\draw [fill](0.15,-0.1) circle [radius=0.11];
\draw [densely dotted]  (0.15,0.9)-- (0.15,-0.1)[postaction={decorate, decoration={markings,mark=at position .55 with {\arrow[black]{stealth}}}}];

\draw  plot[smooth, tension=.7] coordinates {(4,-2) (4.5,-3) (4,-4.5) (3,-5.5)}[postaction={decorate, decoration={markings,mark=at position .5 with {\arrow[black]{stealth}}}}];
\node at (5.5,-4.5) {};
\node at (5.5,-3.5) {};
\draw [fill](5.5,-5) circle [radius=0.11];
\draw [fill](5.5,-4) circle [radius=0.11];
\draw [densely dotted]  (5.5,-4)-- (5.5,-5)[postaction={decorate, decoration={markings,mark=at position .55 with {\arrow[black]{stealth}}}}];

\end{tikzpicture}
\end{matrix}
\end{matrix}
$$
\caption{Change isolated vertices into isolated virtual edges.}
\label{change}
\end{figure}
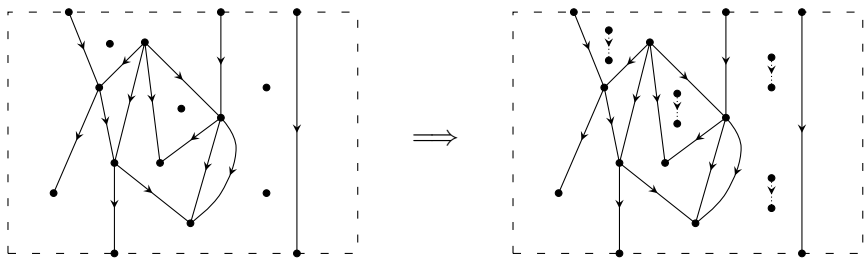

To tackle with appearance of isolated vertices, we can change them into isolated virtual edges, see Fig \ref{change} for an example.

\section{Deformation of upward planar drawing}
Now let us recall Joyal and Street's definition of a \textbf{deformation} of progressive plane graphs.
\begin{defn}[\cite{[JS91]}]\label{def1}
A deformation of two progressive plane graphs $\Gamma_1$, $\Gamma_2$ is a planar isotopy connecting $\Gamma_1$ and $\Gamma_2$ such that each intermediate plane graph is a progressive plane graph.
\end{defn}

There are sufficient algebraic reasons for this definition, which lies in the axioms of monoidal categories and the formalism of string diagrams (we will not give a detailed exposition here, but refer to \cite{[JS91], [S11], [HLY16],[BS10]} and Sect.2 of \cite{[W13]} for readers). The crucial property of this notion --- each intermediate plane graph is a progressive plane graph --- is called \textbf{recumbent}.

Easily, we can generalize this notion for upward plane graphs.

\begin{defn}\label{def2}
A \textbf{deformation} of two upward plane graphs $G_1$, $G_2$ is a planar isotopy connecting $G_1$ and $G_2$ such that each intermediate plane graph is an upward plane graph.
\end{defn}

In other words, we say that a deformation of two upward plane graphs is a recumbent planar isotopy connecting them. Similarly, we have the following notion, which is called an \textbf{upward planar morph} in \cite{[LBFPR18]}.

\begin{defn}\label{def2}
A \textbf{deformation} of two upward planar drawings is a  planar isotopy connecting them such that each intermediate planar drawing is upward.
\end{defn}

Clearly, these notions of deformations define equivalence relations for progressive plane graphs, upward plane graphs and upward planar drawings, and these equivalence relations have strong algebraic motivations and provide a partial answer to $(Q1)$, only from an algebraic point of view.

These deformations are not plain planar isotopies, which defines equivalence relation for plane graphs,  but special planar isotopies with a strong constraint --- the recumbency property.  Then there is a question here: $(Q3)$ in what sense can we say that upward planarity is a topological property? For a complete answer to $(Q1)$, we should answer $(Q3)$, which means to justify the above notions of deformations from a topological view of point.

\section{Bimodal rotation and polarization structure}

A basic result in  topological graph theory \cite{[A96],[MT01]} is the rotation principle (due to Heffter, Edmonds and Ringel), which says that the equivalence class of  an planar drawing of a graph on the plane or an surface can be characterized by a rotation system on this graph. A \textbf{rotation}  of a vertex $v$ of a graph is a cyclic order $\curvearrowright$ on the set $E(v)$ of its incident edges, denoted as $\sigma_v$. A \textbf{rotation system} $\sigma$ of graph $G$ is a collection of rotations, one for each vertex, and we write $\sigma=\{\sigma_v|v\in V(G)\}$.

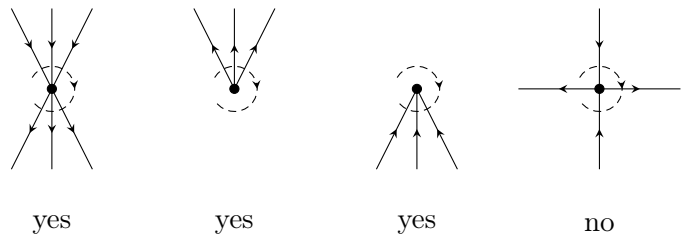
\begin{figure}[H]
\centering
\begin{tikzpicture}[scale=1.2]
\node (v2) at (-2.5,0.5) {};
\node (v8) at (-0.5,0.5) {};
\node (v13) at (1.5,0.5) {};
\node (v17) at (3.5,0.5) {};
\draw [densely dashed][domain=150:-150] plot ({0.25*cos(\x)-2.5}, {0.25*sin(\x)+0.5})[postaction={decorate, decoration={markings,mark=at position .5 with {\arrow[black]{stealth}}}}];
\draw [densely dashed][domain=150:-150] plot ({0.25*cos(\x)-0.5}, {0.25*sin(\x)+0.5})[postaction={decorate, decoration={markings,mark=at position .5 with {\arrow[black]{stealth}}}}];
\draw [densely dashed][domain=150:-150] plot ({0.25*cos(\x)+1.5}, {0.25*sin(\x)+0.5})[postaction={decorate, decoration={markings,mark=at position .5 with {\arrow[black]{stealth}}}}];
\draw [densely dashed][domain=150:-150] plot ({0.25*cos(\x)+3.5}, {0.25*sin(\x)+0.5})[postaction={decorate, decoration={markings,mark=at position .5 with {\arrow[black]{stealth}}}}];
\draw[fill] (v2) circle [radius=0.05];
\draw[fill] (v8) circle [radius=0.05];
\draw[fill] (v13) circle [radius=0.05];
\draw[fill] (v17) circle [radius=0.05];

\node (v1) at (-3,1.5) {};
\node (v3) at (-2.5,1.5) {};
\node (v4) at (-2,1.5) {};
\node (v5) at (-3,-0.5) {};
\node (v6) at (-2.5,-0.5) {};
\node (v7) at (-2,-0.5) {};
\draw  (v1) -- (-2.5,0.5)[postaction={decorate, decoration={markings,mark=at position .5 with {\arrow[black]{stealth}}}}];
\draw  (v3) -- (-2.5,0.5)[postaction={decorate, decoration={markings,mark=at position .5 with {\arrow[black]{stealth}}}}];
\draw  (v4) -- (-2.5,0.5)[postaction={decorate, decoration={markings,mark=at position .5 with {\arrow[black]{stealth}}}}];
\draw (-2.5,0.5) -- (v5)[postaction={decorate, decoration={markings,mark=at position .55 with {\arrow[black]{stealth}}}}];
\draw  (-2.5,0.5) -- (v6)[postaction={decorate, decoration={markings,mark=at position .55 with {\arrow[black]{stealth}}}}];
\draw (-2.5,0.5) -- (v7)[postaction={decorate, decoration={markings,mark=at position .55 with {\arrow[black]{stealth}}}}];
\node (v9) at (-1,1.5) {};
\node (v10) at (0,1.5) {};
\node (v11) at (-0.5,1.5) {};
\draw  (-0.5,0.5) -- (v9)[postaction={decorate, decoration={markings,mark=at position .55 with {\arrow[black]{stealth}}}}];
\draw  (-0.5,0.5) -- (v10)[postaction={decorate, decoration={markings,mark=at position .55 with {\arrow[black]{stealth}}}}];
\draw  (-0.5,0.5) -- (v11)[postaction={decorate, decoration={markings,mark=at position .55 with {\arrow[black]{stealth}}}}];
\node (v12) at (1,-0.5) {};
\node (v14) at (1.5,-0.5) {};
\node (v15) at (2,-0.5) {};
\draw  (v12) -- (1.5,0.5)[postaction={decorate, decoration={markings,mark=at position .5 with {\arrow[black]{stealth}}}}];
\draw  (v14) -- (1.5,0.5)[postaction={decorate, decoration={markings,mark=at position .5 with {\arrow[black]{stealth}}}}];
\draw  (v15) -- (1.5,0.5)[postaction={decorate, decoration={markings,mark=at position .5 with {\arrow[black]{stealth}}}}];
\node (v16) at (3.5,1.5) {};
\node (v20) at (3.5,-0.5) {};
\node (v18) at (2.5,0.5) {};
\node (v19) at (4.5,0.5) {};
\draw  (v16) -- (3.5,0.5)[postaction={decorate, decoration={markings,mark=at position .5 with {\arrow[black]{stealth}}}}];
\draw  (v18) -- (3.5,0.5)[postaction={decorate, decoration={markings,mark=at position .5 with {\arrowreversed[black]{stealth}}}}];
\draw (v19)-- (3.5,0.5)[postaction={decorate, decoration={markings,mark=at position .5 with {\arrowreversed[black]{stealth}}}}];
\draw (v20) -- (3.5,0.5)[postaction={decorate, decoration={markings,mark=at position .5 with {\arrow[black]{stealth}}}}];

\node at (-2.5,-1) {yes};
\node at (-0.5,-1) {yes};
\node at (1.5,-1) {yes};
\node at (3.5,-1) {no};

\end{tikzpicture}

\caption{Bimodal rotations and a non-bimodal rotation }
\label{bimodal}
\end{figure}

For directed graphs, a similar notion should be \textbf{bimodal rotation system}.  For a vertex $v$ of a directed graph, a natural notion is that of a \textbf{bimodal rotation} \cite{[GT95]}, which is a rotation of $v$ such that both the set $I(v)$ of incoming edges and the set $O(v)$ of outgoing edges are intervals with respect to the cyclic order, see Fig \ref{bimodal} for examples and a non-example. Note that for a source or a sink, any of its rotations is naturally bimodal.

\begin{prop}[\cite{[TT86]}]
 Any upward planar graph induces a unique bimodal rotation system on its underlying acyclic directed graph.
\end{prop}

\begin{figure}[htbp]
\centering
$$
\begin{matrix}
\begin{matrix}
\begin{tikzpicture}[scale=0.4]
\node (v2) at (0,0.5) {};
\draw[fill] (v2) circle [radius=0.13];
\node [above] (v1) at (-1.5,2.5) {};
\node [above](v3) at (-0.5,2.5) {};
\node [scale=0.8, below]at (0.5,2.5) {$\cdots$};
\node [above](v4) at (1.5,2.5) {};
\node [below](v5) at (-1.5,-1.5) {};
\node [below](v6) at (-0.5,-1.5) {};
\node[scale=.8, above](v7) at (0.5,-1.5) {$\cdots$};
\node [below](v8) at (1.5,-1.5) {};
\draw  (-1.5,2.5) -- (0,0.5)[postaction={decorate, decoration={markings,mark=at position .5 with {\arrow[black]{stealth}}}}];
\draw  (-0.5,2.5) -- (0,0.5)[postaction={decorate, decoration={markings,mark=at position .5 with {\arrow[black]{stealth}}}}];
\draw  (1.5,2.5) -- (0,0.5)[postaction={decorate, decoration={markings,mark=at position .5 with {\arrow[black]{stealth}}}}];
\draw  (0,0.5) -- (-1.5,-1.5)[postaction={decorate, decoration={markings,mark=at position .55 with {\arrow[black]{stealth}}}}];
\draw  (0,0.5) -- (-0.5,-1.5)[postaction={decorate, decoration={markings,mark=at position .55 with {\arrow[black]{stealth}}}}];
\draw  (0,0.5) -- (1.5,-1.5)[postaction={decorate, decoration={markings,mark=at position .55 with {\arrow[black]{stealth}}}}];
\node[scale=0.8] at (-1.6,1.9) {$1$};
\node [scale=0.8]at (-0.8,2.1) {$2$};
\node [scale=0.8]at (1.5,1.9) {$m$};
\node [scale=0.8]at (-1.5,-0.8) {$1$};
\node [scale=0.8]at (-0.7,-1) {$2$};
\node [scale=0.8]at (1.4,-0.9) {$n$};
\end{tikzpicture}
\end{matrix}&&&&\begin{matrix}
\begin{tikzpicture}[scale=0.4]
\node (v2) at (0,0.5) {};
\draw[fill] (v2) circle [radius=0.13];
\node [above] (v1) at (-1.5,2.5) {};
\node [above](v3) at (-0.5,2.5) {};
\node [scale=0.8, below]at (0.4,2.7) {$\cdots$};
\node [above](v4) at (1.5,2.5) {};
\draw  (-1.68,2.5) -- (0,0.5)[postaction={decorate, decoration={markings,mark=at position .5 with {\arrow[black]{stealth}}}}];
\draw  (-0.75,2.5) -- (0,0.5)[postaction={decorate, decoration={markings,mark=at position .5 with {\arrow[black]{stealth}}}}];
\draw  (1.5,2.5) -- (0,0.5)[postaction={decorate, decoration={markings,mark=at position .5 with {\arrow[black]{stealth}}}}];
\node[scale=0.8] at (-1.6,2.9) {$1$};
\node [scale=0.8]at (-0.8,2.9) {$2$};
\node [scale=0.8]at (1.5,2.9) {$k$};
\end{tikzpicture}
\end{matrix}
&&&&\begin{matrix}
\begin{tikzpicture}[scale=0.4]
\node (v2) at (0,0.5) {};
\draw[fill] (v2) circle [radius=0.13];
\node [above](v4) at (1.5,2.5) {};
\node [below](v5) at (-1.5,-1.5) {};
\node [below](v6) at (-0.5,-1.5) {};
\node[scale=.8, above](v7) at (0.6,-1.7) {$\cdots$};
\node [below](v8) at (1.5,-1.5) {};
\draw  (0,0.5) -- (-1.5,-1.5)[postaction={decorate, decoration={markings,mark=at position .55 with {\arrow[black]{stealth}}}}];
\draw  (0,0.5) -- (-0.5,-1.5)[postaction={decorate, decoration={markings,mark=at position .55 with {\arrow[black]{stealth}}}}];
\draw  (0,0.5) -- (1.8,-1.5)[postaction={decorate, decoration={markings,mark=at position .55 with {\arrow[black]{stealth}}}}];

\node [scale=0.8]at (-1.5,-1.9) {$1$};
\node [scale=0.8]at (-0.6,-1.9) {$2$};
\node [scale=0.8]at (1.8,-1.9) {$l$};
\end{tikzpicture}
\end{matrix}
\end{matrix}
$$
\caption{Examples of polarizations. }
\label{pols}
\end{figure}
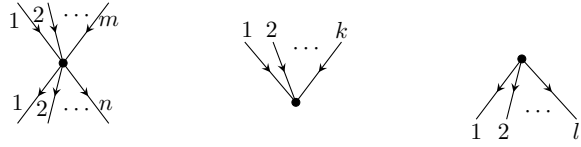

In the theory of graphical calculi, the similar notion is that of a \textbf{polarization structure}. A \textbf{polarization} \cite{[JS91]} of a vertex $v$ of a directed graph is a choice of a linear order on the set $I(v)$ of incoming edges and a linear order on the set $O(v)$ of outgoing edges, see Fig \ref{pols} for an example.  A polarization structure of a directed graph is a collection of polarizations of all its vertices. The notion of a polarization is same as that of \textbf{$2$-linear lists of edges} in \cite{[BB91],[BBLM94]}.

\begin{figure}[H]
\centering
\begin{tikzpicture}[scale=0.9]

\node [right][scale=0.7] at (4.15,-1.5) {$v$};
\draw[fill] (3.5,-1.5) circle [radius=0.055];
\draw [densely dashed][domain=150:-150] plot ({0.55*cos(\x)+3.5}, {0.55*sin(\x)-1.5})[postaction={decorate, decoration={markings,mark=at position .5 with {\arrow[black]{stealth}}}}];
\node [scale=0.5](v17) at (2.5,-0.6) {$1$};
\node [scale=0.5](v19) at (3.2,-0.6) {$2$};
\node [scale=0.5](v20) at (4.5,-0.6) {$k$};
\node[scale=0.5] (v21) at (3.2,-2.4) {$2$};
\node[scale=0.5] (v24) at (2.5,-2.4) {$1$};
\node [scale=0.5](v23) at (4.5,-2.4) {$l$};

\node [scale=0.5]at (3.9,-0.6) {$\cdots$};
\node [scale=0.5](v22) at (3.9,-2.4) {$\cdots$};

\draw  (v17) -- (3.5,-1.5)[postaction={decorate, decoration={markings,mark=at position .4 with {\arrow[black]{stealth}}}}];
\draw  (v19)--(3.5,-1.5)[postaction={decorate, decoration={markings,mark=at position .5 with {\arrow[black]{stealth}}}}];
\draw  (v20) -- (3.5,-1.5)[postaction={decorate, decoration={markings,mark=at position .4 with {\arrow[black]{stealth}}}}];
\draw  (3.5,-1.5) -- (v21)[postaction={decorate, decoration={markings,mark=at position .75 with {\arrow[black]{stealth}}}}];
\draw (3.5,-1.5) -- (v24)[postaction={decorate, decoration={markings,mark=at position .79 with {\arrow[black]{stealth}}}}];
\draw  (3.5,-1.5) -- (v23)[postaction={decorate, decoration={markings,mark=at position .75 with {\arrow[black]{stealth}}}}];

\node[scale=0.7] at (3.5,-0.2) {$I(v)$};
\node[scale=0.7] at (3.5,-2.8) {$O(v)$};

\end{tikzpicture}
\caption{A polarization and its equivalent bimodal rotation of a processive vertex.}
\label{pol}
\end{figure}
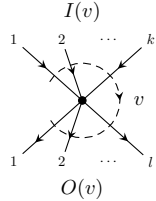

Bimodal rotation has close relations with polarization. As shown in Fig \ref{pol}, any polarization of $v$ induces a bimodal rotation $\curvearrowright$ in the way that $(1)$ for any $f\in O(v)$, $f\curvearrowright h_1\curvearrowright h_2$ $\Longleftrightarrow$ $h_1<h_2$ in $I(v)$; $(2)$ for any $h\in I(v)$, $h\curvearrowright f_1\curvearrowright f_2$ $\Longleftrightarrow$ $f_2<f_1$ in $O(v)$. Conversely, if  $v$ is neither a source nor a sink, then any bimodal rotation of $v$ can define a polarization of $v$ in the above way. However, for a source or a sink, there is no canonical way to define a polarization from a (bimodal) rotation, and the reason is that to break a cyclic order into a linear order we have to choice an element as the first (or last) element and in this case there is no such a natural choice. 

\begin{defn}[\cite{[HLY16]}]
 A vertex of a directed graph is called \textbf{processive} if it is neither a source and a sink.
\end{defn}

The following lemma is clear.
\begin{lem}\label{bp}
For a processive vertices or a leaf, its bimodal rotations and polarizations are in bijective with each other.
\end{lem}

 Just as the case of progressive plane graphs, we have the following result.

\begin{prop}[\cite{[JS91]}]\label{polar}
Any upward plane graph induces a unique polarization structure on the underlying acyclic directed graph.
\end{prop}

\section{NPP-extension and main result}

First, let us recall a notion. 
\begin{defn}[\cite{[HLY16]}]
A \textbf{processive graph} is a non-empty acyclic directed graph with all sources and sinks being leaves (vertices of degree one). 
\end{defn}
A processive graph is a non-empty acyclic directed graph with non-leaves being processive.

By Lemma \ref{bp}, the bimodal rotation systems of a progressive graph are in bijective with its polarization structures .

\begin{figure}[htbp]
\centering
$$
\begin{matrix}
\begin{matrix}
\begin{tikzpicture}[scale=0.4]
\node (v2) at (0,0.5) {};
\draw[fill] (v2) circle [radius=0.13];
\node [left]at (0,0.5) {$s$};

\node (v7) at (0,3) {};
\draw[fill] (v7) circle [radius=0.13];
\node [left]at (0,3) {$l_s$};

\draw [dashed] (0,3) -- (0,0.5)[postaction={decorate, decoration={markings,mark=at position .55 with {\arrow[black]{stealth}}}}];
\node [right]at (0,1.8) {$e_s$};

\node [above](v4) at (1.5,2.5) {};
\node [below](v5) at (-1.5,-1.5) {};
\node [below](v6) at (-0.5,-1.5) {};
\node[scale=.8, above](v7) at (0.6,-1.7) {$\cdots$};
\node [below](v8) at (1.5,-1.5) {};
\draw  (0,0.5) -- (-1.5,-1.5)[postaction={decorate, decoration={markings,mark=at position .55 with {\arrow[black]{stealth}}}}];
\draw  (0,0.5) -- (-0.5,-1.5)[postaction={decorate, decoration={markings,mark=at position .55 with {\arrow[black]{stealth}}}}];
\draw  (0,0.5) -- (1.8,-1.5)[postaction={decorate, decoration={markings,mark=at position .55 with {\arrow[black]{stealth}}}}];

\node [scale=0.8]at (-1.5,-1.9) {};
\node [scale=0.8]at (-0.6,-1.9) {};
\node [scale=0.8]at (1.8,-1.9) {};
\end{tikzpicture}
\end{matrix}&&&&&&
\begin{matrix}
\begin{tikzpicture}[scale=0.4]
\node (v2) at (0,0.5) {};
\draw[fill] (v2) circle [radius=0.13];
\node[left] at (0,0.5) {$t$};

\node (v7) at (0,-2) {};
\draw[fill] (v7) circle [radius=0.13];
\node[left] at (0,-2){$l_t$};
\draw [dashed] (0,0.5) -- (0,-2)[postaction={decorate, decoration={markings,mark=at position .55 with {\arrow[black]{stealth}}}}];
\node [right]at (0,-0.75) {$e_t$};

\node [above] (v1) at (-1.5,2.5) {};
\node [above](v3) at (-0.5,2.5) {};
\node [scale=0.8, below]at (0.4,2.7) {$\cdots$};
\node [above](v4) at (1.5,2.5) {};
\draw  (-1.68,2.5) -- (0,0.5)[postaction={decorate, decoration={markings,mark=at position .5 with {\arrow[black]{stealth}}}}];
\draw  (-0.75,2.5) -- (0,0.5)[postaction={decorate, decoration={markings,mark=at position .5 with {\arrow[black]{stealth}}}}];
\draw  (1.5,2.5) -- (0,0.5)[postaction={decorate, decoration={markings,mark=at position .5 with {\arrow[black]{stealth}}}}];
\node[scale=0.8] at (-1.6,2.9) {};
\node [scale=0.8]at (-0.8,2.9) {};
\node [scale=0.8]at (1.5,2.9) {};
\end{tikzpicture}
\end{matrix}
\end{matrix}
$$
\caption{Normal processive-extensions of $s$ and $t$.}
\label{NE}
\end{figure}
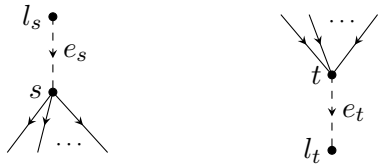

\begin{defn}
For an acyclic directed graph $G$, its \textbf{normal processive-extension} (or \textbf{NP-extension} for short) is the processive graph $\widehat{G}$ obtained from $G$ by the following two kinds of procedures  (see Fig. \ref{NE}):

$(1)$ for each source $s$ of $G$, add one leaf $l_s$ and one directed edge $e_s$ such that $O(l_s)=\{e_s\}=I(s)$;

$(2)$ for each sink $t$ of $G$, add one leaf $l_t$ and one directed edge $e_t$ such that $O(t)=\{e_t\}=I(l_t)$.
\end{defn}
For any acyclic directed graph, its NP-extension exists and is unique, which is just a processive graph obtained by adding an input edge for each source and adding an output edge for each sink.

 For upward plane graphs, we have a similar notion.
\begin{defn}
Let $\Gamma$ be an upward plane graph, its \textbf{normal processive planar extension} (or \textbf{NPP-extension}) is the upward plane graph  $\widehat{\Gamma}$ obtained from $\Gamma$  by adding an input edge for each source and adding an output edge for each sink.
\end{defn}

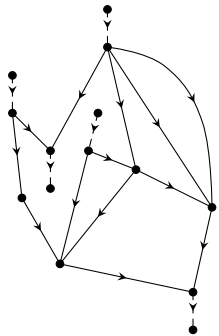
\begin{figure}[H]
\centering
  \begin{tikzpicture}[scale=0.25]

\node (v2) at (-3,-1) {};
\node (v1) at (-2.5,-5.5) {};
\node (v3) at (-1,-3) {};
\node (v4) at (2,2.5) {};
\node (v8) at (1,-3) {};
\node (v5) at (-0.5,-9) {};
\node (v6) at (3.5,-4) {};
\node (v7) at (7.5,-6) {};
\node (v9) at (6.5,-10.5) {};

\draw[fill] (v1) circle [radius=0.2];
\draw[fill] (v2) circle [radius=0.2];
\draw[fill] (v3) circle [radius=0.2];
\draw[fill] (v4) circle [radius=0.2];
\draw[fill] (v5) circle [radius=0.2];
\draw[fill] (v6) circle [radius=0.2];
\draw[fill] (v7) circle [radius=0.2];
\draw[fill] (v8) circle [radius=0.2];
\draw[fill] (v9) circle [radius=0.2];

\draw  (6.5,-10.5) -- (-0.5,-9)[postaction={decorate, decoration={markings,mark=at position .5 with {\arrowreversed[black]{stealth}}}}];
\draw  (-0.5,-9) -- (-2.5,-5.5)[postaction={decorate, decoration={markings,mark=at position .5 with {\arrowreversed[black]{stealth}}}}];
\draw  (-2.5,-5.5) -- (-3,-1)[postaction={decorate, decoration={markings,mark=at position .5 with {\arrowreversed[black]{stealth}}}}];
\draw  (-1,-3) -- (-3,-1)[postaction={decorate, decoration={markings,mark=at position .5 with {\arrowreversed[black]{stealth}}}}];
\draw  (-1,-3) -- (2,2.5)[postaction={decorate, decoration={markings,mark=at position .5 with {\arrowreversed[black]{stealth}}}}];
\draw  (-0.5,-9) -- (1,-3)[postaction={decorate, decoration={markings,mark=at position .5 with {\arrowreversed[black]{stealth}}}}];
\draw  (3.5,-4) -- (1,-3)[postaction={decorate, decoration={markings,mark=at position .5 with {\arrowreversed[black]{stealth}}}}];
\draw  (6.5,-10.5) -- (7.5,-6)[postaction={decorate, decoration={markings,mark=at position .5 with {\arrowreversed[black]{stealth}}}}];
\draw  (-0.5,-9) -- (3.5,-4)[postaction={decorate, decoration={markings,mark=at position .5 with {\arrowreversed[black]{stealth}}}}];
\draw  (3.5,-4) -- (2,2.5)[postaction={decorate, decoration={markings,mark=at position .5 with {\arrowreversed[black]{stealth}}}}];
\draw  (7.5,-6) -- (3.5,-4)[postaction={decorate, decoration={markings,mark=at position .5 with {\arrowreversed[black]{stealth}}}}];
\draw  (7.5,-6) -- (2,2.5)[postaction={decorate, decoration={markings,mark=at position .5 with {\arrowreversed[black]{stealth}}}}];

\node (v10) at (2,2.5) {};
\node (v11) at (7.5,-6) {};
\draw  plot[smooth, tension=.7] coordinates {(v10) (5,1.5) (7,-1.5) (v11)}[postaction={decorate, decoration={markings,mark=at position .5 with {\arrow[black]{stealth}}}}];

\node (v12) at (-3,1) {};
\node (v13) at (2,4.5) {};
\node (v14) at (1.5,-1) {};
\node (v15) at (6.5,-12.5) {};
\draw [densely dashed] (-3,1) -- (-3,-1)[postaction={decorate, decoration={markings,mark=at position .5 with {\arrow[black]{stealth}}}}];
\draw  [densely dashed](2,4.5) -- (2,2.5)[postaction={decorate, decoration={markings,mark=at position .5 with {\arrow[black]{stealth}}}}];
\draw  [densely dashed](1.5,-1) -- (1,-3)[postaction={decorate, decoration={markings,mark=at position .5 with {\arrow[black]{stealth}}}}];
\draw  [densely dashed](6.5,-10.5) -- (6.5,-12.5)[postaction={decorate, decoration={markings,mark=at position .5 with {\arrow[black]{stealth}}}}];
\node (v16) at (-1,-5) {};
\draw  [densely dashed](-1,-3) -- (-1,-5)[postaction={decorate, decoration={markings,mark=at position .5 with {\arrow[black]{stealth}}}}];
\draw[fill] (v12) circle [radius=0.2];
\draw[fill] (v13) circle [radius=0.2];
\draw[fill] (v14) circle [radius=0.2];
\draw[fill] (v15) circle [radius=0.2];
\draw[fill] (v16) circle [radius=0.2];
\end{tikzpicture}
\caption{An example of NPP extension/upward consistent assignment}
\label{NPP}
\end{figure}

Fig \ref{NPP} shows an example of a NPP-extension. The notion of a NPP-extension is equivalent to that of a \textbf{upward consistent assignment} \cite{[BB91],[BBLM94]}, which can be explained by the unit convention in graphical calculus for monoidal categories.

\begin{lem}\label{npp}
For an upward plane graph, its NPP-extension is always exists and is unique. 
\end{lem}

By Lemma \ref{bp} ,Lemma \ref{npp} and Proposition \ref{polar} , we can easily get the following theorem.
\begin{thm}\label{t1}
For an upward plane graph, its NPP-extension and the polarization structure of the underlying acyclic directed graph are equivlent.
\end{thm}

The following theorem is a direct consequence of the  combinatorial characterization of upward planar drawings by Bertolazzi, etc \cite{[BB91],[BBLM94]}.
\begin{thm}\label{t2}
Two upward plane graphs are related by a deformation if and only if their NPP-extensions are related by a plane isotopy.
\end{thm}

Combine Theorem \ref{t1} and Theorem \ref{t2}, we get the following main result.

\begin{thm}
Two upward planar graphs are  related by a deformation if  and only if they are connected by a planar isotopy which preserves the orientation and polarization of the underlying acyclic directed graph.
\end{thm}

This theorem implies that two upward planar drawings of an acyclic directed graph $G$ are equivalent (connected by a deformation) if and only if they are connected by a planar isotopy which preserves the orientation and polarization of $G$. This theorem also gives a  positive answer to Selinger's conjectue \cite{[S11]} with the strategy different from that of Delpeuch and Vicary \cite{[DV18]}.

%%%%%%%%%%%%%%%%%%%%%%%%%%%%%%%%%%%%%%%%%%%%%%%%%%%%%%%%%%%

\par

%%%%%%%%%%%%%%%%%%%%%%%%%% Address & Email %%%%%%%%%%%%%%%%%%%%%%%%%%%%%%%%%%%%%%%%

\textbf{Xuexing Lu}\hfill \\  Email: xxlu@uzz.edu.cn

\end{document}